\RequirePackage{fix-cm}
\documentclass[smallextended]{svjour3}       % onecolumn (second format)
\smartqed  % flush right qed marks, e.g. at end of proof
\usepackage{graphicx}
\usepackage{mathrsfs,color}
\usepackage{mathptmx, amsmath, amsfonts}
\newcommand{\A}{\mathcal{A}}
\newcommand{\Ka}{{\kappa}}
\newcommand{\we}{{w}}
\newcommand{\fd}{{\delta}}

\newcommand{\iprod}[1]{\langle#1\rangle}
\begin{document}
\title{A finite difference method for space fractional differential equations with variable diffusivity coefficient
\thanks{The support of the King Fahd University of Petroleum and Minerals (KFUPM) through the  project No. KAUST005 is gratefully acknowledged. Research reported in this publication was supported by research funding from King Abdullah University of Science and Technology (KAUST).}}
\titlerunning{Finite differences for fractional elliptic models}
% At most 5 thanks
%
\author{K. Mustapha \and  K. Furati \and  O. M. Knio \and   O.P. Le Ma\^{\i}tre}

\institute{K. Mustapha  \and
K. Furati \at Department of Mathematics and Statistics, KFUPM,
Dhahran, 31261, Saudi Arabia, \and O. M. Knio \at Computer, Electrical, Mathermatical Sciences and Engineering Division, KAUST, Thuwal 23955, Saudi Arabia, \and O.P. Le Ma\^{\i}tre \at CNRS, LIMSI, Universit\'e Paris-Scalay, Campus Universitaire - BP 133, F-91403 Orsay, France}

\date{Received: \today / Accepted: date}

\maketitle

\begin{abstract} {Anomalous diffusion is a phenomenon that cannot be modeled accurately by second-order diffusion equations, but is better described by fractional diffusion models. The nonlocal nature of the fractional diffusion operators makes substantially more difficult the mathematical analysis of these models and the establishment of suitable numerical schemes. This paper proposes and analyzes the first finite difference method for solving  {\em variable-coefficient} one-dimensional fractional DEs, with two-sided fractional derivatives (FDs). The proposed scheme combines first-order forward and backward Euler methods for approximating the left-sided FD when the right-sided FD is approximated by two consecutive applications of the first-order backward Euler method. Our  scheme reduces to the standard second-order central difference in the absence of FDs. The existence and uniqueness of the numerical solution are proved,  and truncation errors of order $h$ are demonstrated ($h$  denotes the maximum space step size). The numerical tests illustrate the global $O(h)$ accuracy, except for nonsmooth cases which, as expected, have deteriorated convergence rates.
}
{{\bf AMS subject classifications.} 26A33, 35R09, 65M06, 65M15}
\end{abstract}
%\subjclass{%
%26A33, % fractional derivatives and integrals
%45J05, % integro-partial differential equations
%65M12, % stability and convergence of numerical methods
%65M15, % error bounds
%65M60} % finite elements, Raleigh-Ritz and variational methods

\section{Introduction}\label{sec: intro}
This work aims at constructing and analyzing a finite difference scheme for solving one-dimensional two-sided conservative fractional order differential equations with variable coefficient, $\Ka$, of the form:
\begin{equation}\label{eq: FPDE1}
-\partial_x\left(\Ka(x) \partial_x^{\alpha,\theta}  u(x)\right)=f(x),\quad {\rm for}~~x\in \Omega:=(a,b),
\end{equation}
subject to absorbing boundary  conditions $u=0$ on ${\mathbb R}\backslash \Omega$ and so  $u(a)=u(b)=0$. In \cite{DengLiTianZhang2018},  the authors introduced physically reasonable boundary constraints for different  fractional PDEs.

 In \eqref{eq: FPDE1}, $\alpha \in (0,1)$ is the fractional order exponent, $\Ka$ is the generalized diffusivity coefficient satisfying the positivity assumption $c_0\le \Ka(x) \le c_1$ on  $\Omega$ for some positive constants $c_0$ and $c_1$, $\partial_x$ denotes the first-order derivative, and $\partial_x^{\alpha,\theta}$ the two-sided fractional order differential operator defined by

\[%\begin{equation}\label{Balpha}
\partial^{\alpha,\theta}_x \phi:=\theta {_a{\rm D}_x^{\alpha}} \phi+(1-\theta){_x{\rm D}_b^{\alpha}}\phi.
\]%\end{equation}
Here, $0\le \theta\le 1$  is a parameter describing the relative probabilities of particles to travel ahead or behind the mean displacement,
$_a{\rm D}_x^{\alpha}$ and $_x{\rm D}_b^{\alpha}$ are left-sided (LS) and right-sided (RS) Riemann-Liouville fractional derivatives, defined respectively as
\[  _a{\rm D}_x^{\alpha} v(x) := \frac{\partial }{\partial x}\,  _aI_x^{1-\alpha}v(x)=
\frac{\partial }{\partial x}\int_a^x\omega_{1-\alpha}(x-z)v(z)\,dz,\]
and
\[  _x{\rm D}_b^{\alpha} v(x):=\frac{\partial }{\partial x}\,  _xI_b^{1-\alpha}=
\frac{\partial }{\partial x}\int_x^b \omega_{1-\alpha}(z-x)v(z)dz\,.
\]
In the previous expressions, we denoted $_aI_x^{1-\alpha}$ and $_xI_b^{1-\alpha}$ the LS and RS Riemann-Liouville fractional integrals, respectively, with kernel $\omega_{1-\alpha}(x):=\frac{x^{-\alpha}}{\Gamma(1-\alpha)}$. %where $\Gamma$ is the classical gamma function.

%Without loss of generality, we assume in the following that $d_1=d_2=0$, that is, homogeneous Dirichlet boundary conditions.
%If $d_1\neq 0$ or $d_2 \neq 0$, we substitute
%$$u(x)=w(x)+ \tilde u(x), \quad \tilde u(x) := \frac{x-a}{b-a}d_2+\frac{b-x}{b-a}d_1,$$
%in~\eqref{eq: FPDE1}, and solve for $w$ subject to homogeneous Dirichlet boundary conditions, $w(a)=w(b)=0$, and modified source term function
%$\tilde f := f +\partial_x\left(\Ka \partial_x^{\alpha,\theta} \tilde u\right).$ The practical evaluation of $\tilde f$ can follow the approach outlined in~  %\ref{appen: RHS}.

In the limiting case $\alpha = 1$, the fractional derivative $\partial^\alpha_x$ reduces to $\partial_x$ and the  problem~\eqref{eq: FPDE1} reduces to the classical two-point elliptic boundary value problem, where $-\Ka \partial_x u$ is the ordinary diffusion flux from the Fick'€™s law, Fourier's law, or Newtonian constitutive equation.
An implied assumption is that the rate of diffusion at a certain location is independent of the global structure of the diffusing field.
In the last few decades, an increasing number diffusion processes were found to be non-Fickian, and anomalous diffusion has been experimentally documented in many applications of interest~\cite{BensonWheatcraftMeerschaert2000,MetzlerKlafter2000,podlubny1999} (\textit{e.g.}, viscoelastic materials, subsurface flows and plasma physics). In these situations, the mean square displacement grows in time faster (superdiffusion) or slower (subdiffusion) than that in a normal (Gaussian) diffusion process.
This deviation from normal diffusion can be explained by non-Newtonian mechanics and L\'evy processes.
In such phenomena, the anomalous diffusion rate is affected not only by the local conditions (gradient) but also by the global state of the field.
For instance, the time fractional derivative acting on the diffusion term (subdiffusion)~\cite{MetzlerKlafter2000} accommodates the existence of long-range correlations in the particle dynamics. Similarly, space fractional derivatives, which are suitable for the modeling of superdiffusion processes, account for anomalously large particle jumps at a rate inconsistent with the classical Brownian motion model.
At the macroscopic level, these jumps give rise to a spatial fractional diffusion equation~\cite{BensonWheatcraftMeerschaert2000,Chaves1998}:
\begin{equation}\label{eq: FPDE1 diffusion}
\partial_t u -\partial_x(\Ka \partial_x^{\alpha,\theta}  u)=g.
\end{equation}
In most studies, the diffusion coefficient $\Ka$ is assumed to be constant, and the process to be symmetric~\cite{BensonWheatcraftMeerschaert2000,del-Castillo-NegreteCarrerasLynch2003}. In this case,  $\theta = 1/2$, ~\eqref{eq: FPDE1} reduces to the Riesz fractional derivative of order
$1+\alpha$, and  many numerical methods have been proposed for its solution, see for example \cite{CelikDuman2012,DengHesthaven2013,DingLiChen2014,KharazmiZayernouriKarniadakis2016,LiXu2010,LiuYanKhan2017,LynchEtAl2003,Ortigueira2006,Sousa2009,TadjeranMeerschaertScheffler2006,TianZhouDeng2015,ZhangLiuAnh2010,ZhouTianDeng2013}. However, many practical problems require a model with variable diffusion coefficients $\Ka$~\cite{ChechkinKlafterSokolov2003}, and the asymmetric diffusion process seems inherent in some physical systems~\cite{del-Castillo-Negrete2000,SolomonWeeksSwinney1993}.

The model problem~\eqref{eq: FPDE1} is the steady state form of \eqref{eq: FPDE1 diffusion}.
For a constant diffusivity $\Ka$, the operator $\partial_x(\Ka \partial_x^{\alpha,\theta} )$ is a linear combination of the LS and RS fractional derivatives of order $\alpha+1$.
Let $\iprod{\cdot,\cdot}$ be the  $L_2$-inner product over $\Omega$ and $H^\mu_0(\Omega)$, with $\mu > 1/2$, the fractional Sobolev space of order $\mu$ of functions with zero trace on $\partial \Omega$.
For the Galerkin weak formulation of~\eqref{eq: FPDE1}, we seek the solution $u \in H^{1-\beta}_0(\Omega)$, such that
\begin{equation}\label{eq: ibvp weak}
 \A(u,v)=\iprod{f,v},\quad \forall v \in H^{1-\beta}_0(\Omega),~~{\rm with}~~\beta=(1-\alpha)/2,
\end{equation}
where the bilinear form $\A: H^{1-\beta}_0(\Omega)\times H^{1-\beta}_0(\Omega) \to \mathbb{R}$, is defined by
\[
\A(v,w) :=-\Ka [\theta\iprod{_a{\rm D}_x^{1-\beta} v,\,  _x{\rm D}_1^{1-\beta} w}+(1-\theta)\iprod{_x{\rm D}_1^{1-\beta} v,\, _a{\rm D}_x^{1-\beta} w}].
\]
Ervin and Roop~\cite{ErvinRoop2006} investigated the well-posedness of the Galerkin formulation~\eqref{eq: ibvp weak} for constant $\Ka$.
They proved that the bilinear form $\A$ is then coercive and continuous on $H^{1-\beta}_0(\Omega)\times H^{1-\beta}_0(\Omega) \to \mathbb{R}$, and hence, that~\eqref{eq: ibvp weak} has a unique solution $u \in H^{1-\beta}_0(\Omega)$ in this case.
For a rigorous study of the variational formulation of~\eqref{eq: FPDE1} when $\Ka$ is constant and $\theta=1$, we refer to~\cite{JinLazarovPasciakRundell2015}.

Unfortunately, it was shown in~\cite{WangYang2013} that the Galerkin formulation loses coercivity on $H^{1-\beta}_0(\Omega)\times H^{1-\beta}_0(\Omega) \to \mathbb{R}$ in the variable $\kappa$ case and the authors even propose a counterexample in the case $\theta=1$, see~\cite[Lemma 3.2]{WangYang2013}.
As a result, the weak formulation is not an appropriate framework for variable coefficient $\Ka$, as the Galerkin finite element methods might fail to converge~\cite{WangYangZhu2014}.
As an alternative, a Petrov-Galerkin method was investigated in~\cite{WangYangZhu2015} for the case of LS fractional derivatives ($\theta=1$). For the same setting, a finite difference method was proposed and analyzed in~\cite{StynesGracia2015}.

It is worth to mention that extending existing numerical methods from constant to variable diffusivity is not straightforward, if feasible at all, because of the presence fractional order derivatives.
Similarly, the analyses of the generic problem~\eqref{eq: FPDE1} remain scarce due to the mathematical difficulties induced by LS and RL nonlocal operators, that prevent reusing the results of classical elliptic equations.
Therefore, the main motivation of the present work is to approximate the solution of~\eqref{eq: FPDE1} \textit{via} finite difference methods,
for variable diffusivity $\Ka$ and allowing skewness parameter $0\le \theta\le 1$.
Specifically, we consider numerical schemes based on appropriate combinations of first-order backward and forward differences.
For convenience, we first develop and analyze in Section~\ref{sec: LS} a finite difference scheme for~\eqref{eq: FPDE1} with $\theta=1$, that is, we have to deal with the LS fractional derivative only.
Then, in Section~\ref{sec: RS}, the other limiting case $\theta=0$ with RS fractional derivative only is considered.
The contributions of both LS and RL fractional derivatives are subsequently combined in Section~\ref{sec: two-sided}, to derive the generic finite difference scheme for~\eqref{eq: FPDE1}  that reduces to the classical second-order central difference scheme in the limiting case $\alpha=1$. For each case, we prove the existence and uniqueness of the finite difference solution and show $O(h)$ truncation errors for the resulting schemes, ($h$ is the maximum space step size).
We present several numerical experiments in Section~\ref{sec: numerical results} to support our theoretical convergence results in the case of smooth and non-smooth solutions.
Finally, Section~\ref{sec: conclusion} provides concluding remarks and recommendations for future works.
\section{LS fractional derivative}\label{sec: LS}
For the discretization of the problem, we consider a partition of $\Omega$ with $P$ subintervals $I_{1\le n\le P}$ constructed using a sequence of $(P+1)$ points such that $a=x_0<x_1<x_2<\cdots<x_P=b$. Unless stated otherwise, we shall restrict ourselves to the case of uniform partitions with spatial step size $h=x_n-x_{n-1} = (b-a)/P$. We shall denote $x_{n+1/2} = \frac{x_{n}+x_{n+1}}{2}$ the center of interval $I_{n+1}$.
Denoting $v^n := v(x_n)$, we use the symbol $\fd v^n$ to denote the backward difference defined as $\fd v(x) = \fd v^n := v^n - v^{n-1}$ for  $x \in I_n,$ and  the symbol $\overline  \fd v^n:= v^{n+1/2} - v^{n-1/2},$ to denote the  central difference.
%\subsection{Finite difference scheme}

For the case of LS fractional derivative, that is, $\theta=1,$ Equation~\eqref{eq: FPDE1}   reduces to
\begin{equation} \label{eq: left-sided}
-\partial_x\left({\Ka(x)\, {_a{\rm D}_x^{\alpha}}} u\right)(x)=f(x).
\end{equation}
Using first a forward type difference treatment of the operator $\partial_x$, we propose the following approximation: with $\Ka^{n+1/2} := \Ka(x_{n+1/2})$,
\begin{equation}\label{eq: integrated FP1}
\partial_x\left(\Ka\, {_a{\rm D}_x^{\alpha}} u(x_n)\right)\approx h^{-1}\left[
	\Ka^{n+1/2}{_a{\rm D}_x^{\alpha}} u(x_{n+1}) - \Ka^{n-1/2}{_a{\rm D}_x^{\alpha}} u(x_{n}) \right].
\end{equation}
Observe that the proposed scheme involves a half-cell shift in the localization of the values of $\Ka$,  resembling the case of the classical second-order elliptic equation.

Remarking that ${_a{\rm D}_x^{\alpha}} u={_aI_x^{1-\alpha}} u'$, because $u(0)=0$, equation~\eqref{eq: integrated FP1} can be recast as
%\begin{equation}\label{eq: integrated FP}
\[\partial_x(\Ka\, {_a{\rm D}_x^{\alpha}} u)(x_n)\approx h^{-1}[\Ka^{n+1/2}{_aI_x^{1-\alpha}} u'(x_{n+1})-\Ka^{n-1/2}{_aI_x^{1-\alpha}} u'(x_{n})].
\]%\end{equation}
Applying now the backward difference approximation to the derivatives inside the integrals, results in
\[%\begin{equation}\label{eq: space stepping}
\partial_x(\Ka {_a{\rm D}_x^{\alpha}} u)(x_n)\approx h^{-2}[\Ka^{n+1/2}({_aI_x^{1-\alpha}} \fd u)(x_{n+1})-\Ka^{n-1/2}({_aI_x^{1-\alpha}} \fd u)(x_{n})],
\]%\end{equation}
for $n=1,\dots,P-1$.  In addition, we have
\begin{equation}\label{eq: LSintegral approx}
\begin{aligned}
_aI_x^{1-\alpha} \fd u(x_n)
&	=\sum_{j=1}^n\int_{I_j}\omega_{1-\alpha}(x_n-s)\fd u^j\,ds
	=\omega_{2-\alpha}(h)\sum_{j=1}^n\we_{n,j}\fd u^j \\
& =\omega_{2-\alpha}(h)\Big(\sum_{j=1}^{n-1}[\we_{n,j}-\we_{n,j+1}] u^j+u^n\Big),
\end{aligned}
\end{equation}
with the weights defined as
\begin{equation}\label{eq: weights}
\we_{n,j} := (n+1-j)^{1-\alpha}-(n-j)^{1-\alpha}
	\quad\text{for $ n\ge j\ge 1$.}
\end{equation}
We denote by $U^n\approx u^n$ the finite difference solution, which for the model problem in~\eqref{eq: left-sided} is required to satisfy
\begin{equation}\label{eq: scheme LS}
\Ka^{n-1/2}({_aI_x^{1-\alpha}} \fd U)(x_{n})-\Ka^{n+1/2}({_aI_x^{1-\alpha}} \fd U)(x_{n+1})=h^2 f^n,\quad n=1,\cdots,P-1,
\end{equation}
 with $U^0=U^P=0$. Using~\eqref{eq: LSintegral approx}, the finite difference scheme can be recast as
\begin{equation}\label{eq: scheme LS rewrite}
\Ka^{n-1/2}\sum_{j=1}^{n}\we_{n,j}\fd U^j-\Ka^{n+1/2}\sum_{j=1}^{n+1}\we_{n+1,j}\fd U^j=\tilde f_h^n,
\end{equation}
with the modified right-hand-side
\begin{equation}\label{eq: rhsLS}
\tilde f_h^n:=\frac{h^2}{\omega_{2-\alpha}(h)} f^n.
\end{equation}
For computational convenience,  \eqref{eq: scheme LS} can be expressed in a compact form  as
%\begin{multline*}%\label{eq: linear system}
%\sum_{j=1}^{n-1}\Big(\Ka^{n-1/2}[\we_{n,j}-\we_{n,j+1}]-\Ka^{n+1/2}[\we_{n+1,j}-\we_{n+1,j+1}]\Big) U^j\\
%+\Big(\Ka^{n-1/2}-\Ka^{n+1/2}(2^{1-\alpha}-2) \Big)U^n-\Ka^{n+1/2}U^{n+1}=\tilde f_h^n ,
%\end{multline*}
%or in the compact form:
\begin{equation*}%\label{eq: scheme LS-2}
\sum_{j=1}^{n}\Big(a_{n,j}-a_{n+1,j}\Big) U^j-\Ka^{n+1/2} U^{n+1}=\tilde f_h^n,~~{\rm for}~~n=1,\cdots,P-1,\end{equation*}
where $a_{n,n} = \Ka^{n-1/2}$ and  $a_{n,j}=\Ka^{n-1/2}[\we_{n,j}-\we_{n-1,j}]$ for $j<n.$
	
The finite difference solution is then obtained solving the $(P-1)$-by-$(P-1)$ linear system
${\bf B}_L{\bf  U}={\bf F}$, where  ${\bf  U}=[U^1,U^2,\cdots,U^{P-1}]^T$, ${\bf F}=[\tilde f^1_h, \tilde f^2_h,\cdots,\tilde f^{P-1}_h]^T$, and the matrix ${\bf B}_L=[c_{n,j}]$  having lower-triagonal entries
\begin{equation*}
	c_{n,j} = \begin{cases} \Ka^{n-1/2} +\Ka^{n+1/2}[2-2^{1-\alpha}] & j=n, \cr
	a_{n,j} - a_{n+1,j} & j<n,
	\end{cases}
\end{equation*}
while $c_{n,n+1}=-\kappa^{n+1/2}$  and all other entries are zeros. Note that  for the case of a constant diffusivity, the matrix ${\bf B}_L$ reduces to the Toeplitz form.
 \begin{remark}
As mentioned earlier, in the limiting case $\alpha = 1$, equation~\eqref{eq: FPDE1} reduces to  $-\partial_x(\Ka \partial_x  u)=f.$ Furthermore, the finite difference scheme~\eqref{eq: scheme LS} reduces to
\begin{equation*}%\label{eq: scheme LS alpha=1}
\Ka^{n+1/2} \fd U^{n+1}-\Ka^{n-1/2}\fd U^{n}=h^2 f^n,
\end{equation*}
for $n=1,\cdots,P-1$.
This is the classical second order difference scheme for elliptic problems.
In this case, one can easily check that the system matrix ${\bf B}_L$ becomes tridiagonal and symmetric, with entries
$c_{i,j}=0$ for $|i-j|>2$,  $c_{i,i+1}=-\Ka^{i+1/2}$, $c_{i,i}= \Ka^{i-1/2}+\Ka^{i+1/2}$ and $c_{i,i-1} = -\Ka^{i-1/2}$.
\end{remark}

\begin{lemma} \label{sec: existence LS}
For $1\le n\le P,$ the finite difference solution $U^n$ of \eqref{eq: scheme LS} exists and is unique.
\end{lemma}
{\em Proof.} Since the finite difference solution $U^n$ satisfies a square linear system of equations, % ${\bf B}_L {\bf U},={\bf F}$,
the existence of $U^n$ follows  from its uniqueness. To prove uniqueness, we need to show that the finite difference solution is identically zero when $f=0$, that is when the system right-hand-side is zero, that is $f^j=0$ for $j=1,\cdots,P-1$ in \eqref{eq: scheme LS}. To do so, sum \eqref{eq: scheme LS rewrite} over index $n$, leading to
\[
\sum_{n=1}^m\Ka^{n-1/2}\sum_{j=1}^{n}\we_{n,j}\fd U^j-\sum_{n=1}^m\Ka^{n+1/2}\sum_{j=1}^{n+1}\we_{n+1,j}\fd U^j=0,\]
and consequently,
\[
\sum_{n=0}^{m-1}\Ka^{n+1/2}\sum_{j=1}^{n+1}\we_{n+1,j}\fd U^j-\sum_{n=1}^m\Ka^{n+1/2}\sum_{j=1}^{n+1}\we_{n+1,j}\fd U^j=0.\]
After simplifying, we conclude that
\begin{equation}\label{eq: matrix0}
\Ka^{m+1/2}\sum_{j=1}^{m+1}\we_{m+1,j}\fd U^j=\Ka^{1/2}\fd U^1,\quad {\rm for}~~1\le m\le P-1,\end{equation}
which can alternatively be expressed as
\begin{equation}\label{eq: matrix1}
	{\bf W}_\alpha \Phi =\fd U^1{\bf K},
\end{equation}
where  $\Phi =[\fd U^1,\fd U^2,\cdots,\fd U^P]^T$, ${\bf K} =[k_1,k_2,\cdots,k_P]^T$ with $k_j=\Ka^{1/2}/\Ka^{j-1/2}$, and \begin{equation}\label{eq: Walpha}
{\bf W}_\alpha= \begin{bmatrix}
b_0&0&0&0&0&\cdots&0\\
b_1&b_0&0&0&0&\cdots&0\\
b_2&b_1&b_0&0&0&\cdots&0\\
b_3&b_2&b_1&b_0&0&\cdots&0\\
\vdots&\vdots&\vdots&\vdots &\vdots&\cdots&\vdots \\
b_{P-1}&b_{P-2}&w_{P-3}&b_{P-4}&\cdots&b_1&b_0
\end{bmatrix},\end{equation}
with $b_0=1$ and $b_j=(j+1)^{1-\alpha}-j^{1-\alpha}>0$  for $ j\ge 1.$ Since ${\bf W}_\alpha$ is a   nonsingular lower triangular Toeplitz matrix, its inverse, denoted by ${\bf E}_\alpha$, is also a  lower triangular Toeplitz matrix with elements \[e_0=\frac{1}{b_0}=1,\quad{\rm  and}\quad  e_j=-\sum_{i=0}^{j-1}b_{j-i}\,e_i,\quad {\rm for}~~ j\ge 1.\]
Now, from \eqref{eq: matrix1}, $\Phi ={\bf E}_\alpha{\bf K}\fd U^1$ and thus   $\fd U^j=\fd U^1\sum_{i=1}^{j} e_{j-i}\, k_i\,.$ Since $\sum_{j=1}^P\fd U^j=0$ (because $U^0=U^P=0$), \begin{equation}\label{eq: SDM}
\fd U^1\sum_{j=1}^P\sum_{i=1}^{j} e_{j-i}\, k_i=\fd U^1\sum_{i=1}^Pk_i\sum_{j=i}^{P} e_{j-i} =\fd U^1\sum_{i=1}^Pk_i\sum_{j=0}^{P-i} e_{j} =0.\end{equation}
On the other hand, the sequence $\{b_j\}_{j\ge 0}$  is positive, slowly decaying ($\lim_{j\rightarrow\infty}b_j=0$ and $\sum_{j=1}^\infty |b_k|=\infty$) and is   strictly {\em log-convex} ($b_j^2<b_{j-1}b_{j+1}$ for $j\ge 1$). Then, we deduce that $e_n<0$ and  $\sum_{j=0}^n e_j> 0$ for $n\ge 1,$ see \cite[Theorem 22]{Hardy1949} or \cite[Theorem 2.2 and Lemma 2.4]{FordSavostyanovZamarashkin2014}.  Using this in \eqref{eq: SDM} and also using the fact that $k_i>0$ for $i\ge 1,$ yield  $\fd U^1=0.$ Therefore, by \eqref{eq: matrix1}, $\Phi \equiv {\bf 0}$  (${\bf W}_\alpha$ is nonsingular). Consequently, the finite difference solution $U^n$ is identically zero, for $1\le n\le P-1$, because $U^0=U^P=0$.
This completes the proof of the uniqueness of the numerical solution $U$.$\quad\Box$

We now turn to establishing the truncation error of the proposed scheme.
From~\eqref{eq: left-sided} and~\eqref{eq: scheme LS}, the truncation error $T_h^n$ is given by $T_h^n=\partial_x(\Ka {_a{\rm D}_x^{\alpha}} u)(x_n)-Q^n_h,$
where
\[ Q^n_h
	=\frac{1}{h^2}\Big[\Ka^{n+1/2}{_aI_x^{1-\alpha}} \fd u(x_{n+1})-\Ka^{n-1/2}{_aI_x^{1-\alpha}} \fd u(x_{n})\Big].
\]
Since $\partial_x(\Ka {_a{\rm D}_x^{\alpha}} u)(x_n)= [f(x_n)-f(x)]+\partial_x(\Ka {_a{\rm D}_x^{\alpha}} u)(x),$
\[
\int_{I_{n+1}} \partial_x(\Ka {_a{\rm D}_x^{\alpha}} u)(x_n)\,dx
=- \frac{h^2}{2} f'(\zeta_n) + \Ka^{n+1} {_aI_x^{1-\alpha}} u'(x_{n+1})-\Ka^{n} {_aI_x^{1-\alpha}} u'(x_{n}),\]
for some $\zeta_n\in I_{n+1}$, and thus,
\[
	T_h^n=- \frac{h}{2} f'(\zeta_n) +\frac{1}{h}\Big[\Ka^{n+1} {_aI_x^{1-\alpha}} u'(x_{n+1})-\Ka^{n} {_aI_x^{1-\alpha}} u'(x_{n})\Big]-Q^n_h.
\]
\begin{theorem} \label{sec: truncation LS}
Assume that $f\in C^1(\overline \Omega)$, $\Ka\in C^2(\overline \Omega)$ and $u\in C^3(\overline \Omega)$. Then
\[
	T_h^n = O(h)(1+{(x_n-a)}^{-\alpha}), \quad \mbox{for }  1 \le n\le P-1.
\]
That is,  the truncation error is of order $h$ for $x_n$ not too close to the left boundary.\end{theorem}
{\em Proof.} Using the change of variable $s=q+h$, we observe that
\begin{align*}{_aI_x^{1-\alpha}} u'(x_{n+1})&=
	\sum_{j=1}^{n+1}\int_{I_j}\omega_{1-\alpha}(x_{n+1}-s)u'(s)\,ds\\
&=\int_{I_1}\omega_{1-\alpha}(x_{n+1}-s)u'(s)\,ds+\sum_{j=1}^{n}\int_{I_j}\omega_{1-\alpha}(x_{n}-q)u'(q+h)\,dq.
\end{align*}
Similarly, for the backward difference we have
\begin{align*}
	{_aI_x^{1-\alpha}} \fd u(x_{n+1})&=
	\sum_{j=1}^{n+1}\int_{I_j}\omega_{1-\alpha}(x_{n+1}-s)\fd u^j\,ds\\
&=\fd u^{1}\int_{I_1}\omega_{1-\alpha}(x_{n+1}-s)\,ds+\sum_{j=1}^{n}\fd u^{j+1}\int_{I_j}\omega_{1-\alpha}(x_{n}-q)\,dq.
\end{align*}
Therefore, the truncation error can be rewritten as
\[
T_h^n=- \frac{h}{2} f'(\zeta_n)+E_1^n+\sum_{j=1}^{n}\int_{I_j}\omega_{1-\alpha}(x_{n}-q) E_2^{n,j}(q)\,dq,\quad {\rm for}~~n\ge 1,
\]
where
\[
	E_1^n := h^{-1}\int_{I_1}\omega_{1-\alpha}(x_{n+1}-s)[\Ka^{n+1}u'(s)-h^{-1}\Ka^{n+1/2}u^1]ds,
\]
and
\[
E_2^{n,j}(q):=\frac{\Ka^{n+1}u'(q+h)-\Ka^n u'(q)}{h}-\frac{\Ka^{n+1/2}\fd u^{j+1}-\Ka^{n-1/2}\fd u^{j}}{h^2}\,.
\]
Focusing on the second error contribution, $E^n_1$, we observe that for sufficient smoothness, specifically for $\Ka \in C^1(I_{n+1})$ and $u\in C^2(a,x_1]$, we have  (at leading order)
\begin{align*}
	\Ka^{n+1}u'(s)-h^{-1}\Ka^{n+1/2}u^1&=\Ka^{n+1}u'(s)-h^{-1}[\Ka^{n+1}+O(h)][h\, u'(x_1)+O(h^2)]\\
&=\Ka^{n+1}[u'(s)-u'(x_1)]+O(h)=O(h).
\end{align*}
Consequently, an application of the mean value theorem for integral yields
\begin{align}
	E^n_1 &= O(1)\int_{I_1} {\omega_{1-\alpha}}(x_{n+1}-s)ds = O(h)(x_{n+1}-\xi)^{-\alpha},\quad {\rm for~ some}~~ \xi \in I_1.
\end{align}
  Regarding the last error contribution in $T^n_h$ above, we first remark that for any $q \in (x_{j-1},x_j)$, one has
\[
	\Ka^{n+1}u'(q+h)-\Ka^n u'(q)= \Ka^n [u'(q+h){-u'(q)}]+\delta \Ka^{n+1}u'(q+h),
\]
and that, for $\Ka\in C^2[x_{n-1},x_{n+1}]$ and $u\in C^3[x_{j-1},x_{j+1}]$, Taylor series expansions give
\begin{align*}
\Ka^{n+1/2}\fd u^{j+1}-\Ka^{n-1/2}\fd u^{j}
&= [(\Ka^{n+1/2}-\Ka^n)+\Ka^n][\fd u^{j+1}-\fd u^{j}]+\overline\delta\Ka^{n}\fd u^{j}\\
&=h^2[\frac{h}{2}\Ka'(x_n)+\Ka^n] u''(x_j)+h^2 \Ka'(x_n)u'(x_j)+O(h^3),
\end{align*}
Gathering the previous results, we obtain for $E^{n,j}_2$
\begin{align*}
E_2^{n,j}&(q)  =  h^{-1}\Ka^{n}[u'(q+h)- u'(q)-hu''(x_j)] \\
&~~\quad  +h^{-1}[\delta\Ka^{n+1}-h\Ka'(x_n)]u'(q+h)
+\Ka'(x_n)[u'(q+h)-u'(x_j)]  \\
&= -h^{-1}\Ka^{n}\int_q^{q+h}\int_t^{x_j} u'''(x)\,dx\,dt   +\frac{h}{2}\Ka''(\xi^n)u'(q+h)
+\Ka'(x_n)\int_{x_j}^{q+h}u''(x)\,dx,
\end{align*}
for some $\xi^n \in I_{n+1}.$ This shows that the first double integral term is $O(h^2)$ when $u\in C^3[x_{j-1},x_{j+1}]$, whereas the second term is $O(h)$ for $\Ka\in C^2(\overline I_{n+1})$ and $u\in C^1 (\overline I_{j+1})$ and the third one is $O(h)$ for $\Ka \in C^1(\overline I_{n+1})$ and $u\in C^2(\overline I_{j+1})$. This  leads to the conclusion that the last error contribution to $T^n_h$ is $O(h)$. Putting all these estimates together, we obtain the desired result. $\quad \Box$

\section{RS fractional derivative}\label{sec: RS}
In this section, we focus on the finite difference approximation of problem~\eqref{eq: FPDE1} when $\theta=0$, that is, the RS fractional elliptic problem:
\begin{equation} \label{eq: right-sided}
-\partial_x({\Ka\, {_x{\rm D}_b^{\alpha}}} u)(x)=f(x).
\end{equation}
We shall rely on the same notations as in the previous section.

%\subsection{Finite difference scheme}
Contrary to the case of the LS fractional derivative, we propose a backward difference type treatment for the differential operator $\partial_x$, and
consider the approximation
\begin{equation*}%\label{eq: differential right FP}
\partial_x(\Ka\, {_x{\rm D}_b^{\alpha}} u)(x_n)\approx h^{-1}[\Ka^{n+1/2}{_x{\rm D}_b^{\alpha}} u(x_{n})-\Ka^{n-1/2}{_x{\rm D}_b^{\alpha}} u(x_{n-1})].
\end{equation*}
Again, observe the shift in the evaluation points for $\Ka$ (at the cell centers) compared to fractional differential operator (at the mesh point),
which is crucial to ensure the recovery of the classical second order scheme when $\alpha\rightarrow 1$.
Noting that ${_x{\rm D}_b^{\alpha}} u={_xI_b^{1-\alpha}} u'$, because $u(1)=0$, we have
\begin{equation*}\label{eq: integrated right FP}
\partial_x(\Ka\, {_x{\rm D}_b^{\alpha}} u)(x_n)\approx h^{-1}\left[\Ka^{n+1/2}{_xI_b^{1-\alpha}} u'(x_{n})-\Ka^{n-1/2}{_xI_b^{1-\alpha}} u'(x_{n-1})\right].
\end{equation*}
Applying the backward difference to the derivatives inside the integrals, one gets
\begin{equation*}%\label{eq: space stepping right}
\partial_x(\Ka\, {_x{\rm D}_b^{\alpha}} u)(x_n)\approx h^{-2}\left[\Ka^{n+1/2}({_xI_b^{1-\alpha}} \fd u)(x_{n})-\Ka^{n-1/2}({_xI_b^{1-\alpha}} \fd u)(x_{n-1})\right],
\end{equation*}
for $n=1,\cdots,P-1.$ The finite difference solution $U^n\approx u^n$ of the (RS) fractional model problem~\eqref{eq: right-sided} satisfies the system:
\begin{equation}\label{eq: scheme RS}
\Ka^{n-1/2}({_xI_b^{1-\alpha}} \fd U)(x_{n-1})-\Ka^{n+1/2}({_xI_b^{1-\alpha}} \fd U)(x_{n})=h^2 f^n,
\end{equation}
for $n=1,\cdots,P-1$, complemented by the boundary conditions $U^0=U^P=0$.

Further, application of the integral form of the RS Riemann-Liouville fractional derivative to the finite difference, $\fd v$, yields:
\begin{equation*}%\label{eq: expansion 1 right}
_xI_b^{1-\alpha} \fd v(x_{n-1})
	=\sum_{j=n}^P\int_{I_j}\omega_{1-\alpha}(s-x_{n-1})\fd v^j\,ds
	=\omega_{2-\alpha}(h)\sum_{j=n}^P\we_{j,n}\fd v^j,
\end{equation*}
such that the numerical scheme~\eqref{eq: scheme RS} can be expressed as
\begin{equation}\label{eq: scheme RS rewrite}
	\Ka^{n-1/2}\sum_{j=n}^{P}\we_{j,n}\fd U^j-\Ka^{n+1/2}\sum_{j=n+1}^{P}\we_{j,n+1}\fd U^j=\tilde f_h^n.
\end{equation}
In \eqref{eq: scheme RS}, the weights $\we_{n,j}$ and modified right-hand side $\tilde f_h^n$ follow the definitions of the previous section, see equations~\eqref{eq: weights} and~\eqref{eq: rhsLS} respectively.
Making use of the equality
\[
	\sum_{j=n}^P\we_{j,n}\fd v^j= \sum_{j=n}^{P-1}[\we_{j,n}-\we_{j+1,n}] v^j-\we_{n,n}v^{n-1},
\]
%equation~\eqref{eq: scheme RS rewrite} becomes
%\begin{multline*}%\label{eq: linear system}
%\sum_{j=n+1}^{P-1}\Big(\Ka^{n-1/2}[\we_{j,n}-\we_{j+1,n}]-\Ka^{n+1/2}[\we_{j,n+1}-\we_{j+1,n+1}]\Big) U^j\\
%+\Big(\Ka^{n-1/2}[\we_{n,n}-\we_{n+1,n}]+\Ka^{n+1/2} \we_{n+1,n+1} \Big)U^n-\Ka^{n-1/2}U^{n-1}=\tilde f_h^n.
%\end{multline*}
%Therefore,
the finite difference scheme~\eqref{eq: scheme RS} can be rewritten as
\begin{equation*}%\label{eq: scheme LS-2}
\sum_{j=n}^{P-1}\Big(b_{jn}-b_{j,n+1}\Big) U^j-\Ka^{n-1/2}U^{n-1}=\tilde f_h^n, \quad n=1,\cdots,P-1,
\end{equation*}
where $b_{n,n+1} =-\Ka^{n+1/2}$ and $b_{j,n}=  \Ka^{n-1/2}[\we_{j,n}-\we_{j,n-1}]$ for $j\ge n.$
%\begin{eqnarray}
%	b_{j,n}& = & \Ka^{n-1/2}[\we_{j,n}-\we_{j,n-1}] \nonumber \\
%	           & = & -\Ka^{n-1/2}[(j+2-n)^{1-\alpha}-2(j+1-n)^{1-\alpha}+(j-n)^{1-\alpha}], \nonumber
%\end{eqnarray}
%for $j\ge n$.

The finite difference solution of the RS fractional diffusion problem is thus obtained by solving the
$(P-1)$-by-$(P-1)$ linear system ${\bf B}_R{\bf  U}={\bf F}$, with the system matrix ${\bf B}_R=[d_{n,j}]$  having upper-triagonal entries
% $d_{n,n}=\Ka^{n+1/2}-\Ka^{n-1/2}[2^{1-\alpha}-2]$ and  $d_{n,j}=b_{j,n}-b_{j,n+1}$ when $j>n$.
\begin{equation*}
d_{n,j} = \begin{cases}
-\Ka^{n-1/2}\we_{j,n-1} + (\Ka^{n-1/2}+\Ka^{n+1/2})\we_{j,n} - \Ka^{n+1/2}\we_{j,n+1}, & j>n, \cr
\Ka^{n+1/2}+\Ka^{n-1/2}[2-2^{1-\alpha}], & j=n,
\end{cases}
\end{equation*}
while $d_{n+1,n}=-\kappa^{n+1/2}$ and all other entries are zeros.
%\subsection{Existence and uniqueness}\label{sec: existence RS}
\begin{lemma} \label{sec: existence RS}
The finite difference solution $U^n$  to the RS scheme \eqref{eq: scheme RS} exists and is unique.
\end{lemma}
{\em Proof.}  As in the case of the LS fractional derivative, the existence of the solution $U^n$ of  \eqref{eq: scheme RS} follows from its uniqueness, and it is sufficient to show that the finite difference solution is identically zero when $f^n=0$ for $n=1,\cdots,P-1$. To do so, we follow the same path as in Lemma \ref{sec: existence LS}. Summing \eqref{eq: scheme RS rewrite} over the index $n$, we get
\[
	\sum_{n=m}^{P-1}\Ka^{n-1/2}\sum_{j=n}^{P}\we_{j,n}\fd U^j-\sum_{n=m}^{P-1}\Ka^{n+1/2}\sum_{j=n+1}^{P}\we_{j,n+1}\fd U^j=0.
\]
The second sum equals  $\sum_{n=m+1}^P\Ka^{n-1/2}\sum_{j=n}^{P}\we_{j,n}\fd U^j$, and so,
\[
	\Ka^{m-1/2}\sum_{j=m}^{P}\we_{j,m}\fd U^j-\Ka^{P-1/2}\fd U^P=0.
\]
and it ensues that
\[
	\Ka^{n-1/2}\sum_{j=n}^{P}\we_{j,n}\fd U^j=\Ka^{P-1/2}\fd U^P,\quad{\rm for}~~1\le n\le P.
\]
This equation can be cast in the matrix form,
 \begin{equation}\label{eq: matrix2}
{\bf W}_\alpha^T \Phi =\fd U^P\, \hat {\bf K} \Longleftrightarrow {\bf W}_\alpha \Phi =\fd U^P\, \tilde {\bf K},
\end{equation}
with the same matrix ${\bf W}_\alpha$ as in equation~\eqref{eq: matrix1}, whereas   $\hat {\bf K} =[\frac{\Ka^{P-1/2}}{\Ka^{1/2}},\frac{\Ka^{P-1/2}}{\Ka^{3/2}},\cdots,1]^T$ and $\tilde {\bf K} =[1,\frac{\Ka^{P-1/2}}{\Ka^{P-3/2}},\cdots,\frac{\Ka^{P-1/2}}{\Ka^{1/2}}]^T$. Since \eqref{eq: matrix1} and \eqref{eq: matrix2} have the same form, by following the derivation in Lemma \ref{sec: existence LS}, we deduce that $\fd U^P=0.$
It is again immediate to conclude  from \eqref{eq: matrix2} that $\Phi \equiv {\bf 0}$ because ${\bf W}_\alpha$ is nonsingular. Consequently, the  finite difference solution $U^n=0$ for $1\le n\le P-1$ because $U^0=U^P=0$.
Therefore, the solution to the RS scheme \eqref{eq: scheme RS} exists and is unique. $\quad \Box$

%\subsection{Truncation error}\label{sec: truncation RS}
Next, we study  the truncation  error $T_h^n$ of the proposed finite difference discretization of problem~\eqref{eq: right-sided}. The truncation error  in this case is
\[
T_h^n=\partial_x(\Ka\, {_x{\rm D}_b^{\alpha}} u)(x_n)- Q_h^n
\]
where
\[
	Q^n_h=\frac{1}{h^2}\Big(\Ka^{n+1/2}({_xI_b^{1-\alpha}} \fd u)(x_{n})-\Ka^{n-1/2}({_xI_b^{1-\alpha}} \fd u)(x_{n-1})\Big)
\]
is the proposed finite difference approximation of the RS operator.
Regarding the continuous part, we proceed with a procedure similar to the LS~case, to get
\begin{align*}
h\,\partial_x(\Ka\, {_x{\rm D}_b^{\alpha}} u)(x_n) =  \int_{I_{n}}\partial_x(\Ka\, {_x{\rm D}_b^{\alpha}} u)(x_n)\,dx
%&= \int_{I_{n}}[f(x_n)-f(x)]\,dx+\int_{I_n}\partial_x(\Ka\, {_x{\rm D}_b^{\alpha}} u)%(x)\,dx\\
&= \frac{h^2}{2} f'(\zeta_n)+hG_h^{n},\quad{\rm for~ some}~~ \zeta_n\in I_n,
\end{align*}
where
\begin{equation*}
G_h^{n} = \frac{1}{h}\Big[\Ka^{n} {_xI_b^{1-\alpha}} u'(x_{n})-\Ka^{n-1} {_xI_b^{1-\alpha}} u'(x_{n-1})\Big].\end{equation*}
Consequently,
\begin{equation}\label{eq: imp estimate}
	T_h^n=O(h)+ G_h^{n}-Q^n_h\,.
\end{equation}

  \begin{theorem} \label{sec: truncation RS}
Assume that $f\in C^1(\overline \Omega)$, $\Ka\in C^2(\overline \Omega)$ and $u\in C^3(\overline \Omega)$. Then
\[
	T_h^n = O(h)(1+{(b-x_{n-1})}^{-\alpha}), \quad \mbox{for }  1 \le n\le P-1.
\]
That is,  the truncation error is of order $h$ for $x_n$ not too close to the right boundary.\end{theorem}
{\rm Proof.} Noting first that
\begin{equation}\label{exp_chgvar1}
\begin{aligned}
	_xI_b^{1-\alpha} &u'(x_{n-1}) =
	\sum_{j=n}^{P}\int_{I_j}\omega_{1-\alpha}(s-x_{n-1})u'(s)\,ds\\
&\quad =  \int_{I_P}\omega_{1-\alpha}(s-x_{n-1})u'(s)\,ds  +\sum_{j=n+1}^{P}\int_{I_j}\omega_{1-\alpha}(q-x_{n})u'(q-h)\,dq,
\end{aligned}
\end{equation}
and
\begin{equation}\label{exp_chgvar2}
\begin{aligned}	_xI_b^{1-\alpha} \fd &u(x_{n-1}) =
	\sum_{j=n}^{P}\int_{I_j}\omega_{1-\alpha}(s-x_{n-1})\fd u^j\,ds\\
 &\quad=  \fd u^{P}\int_{I_P}\omega_{1-\alpha}(s-x_{n-1})\,ds + \sum_{j=n+1}^{P}\fd u^{j-1}\int_{I_j}\omega_{1-\alpha}(q-x_{n})\,dq.
\end{aligned}
\end{equation}
On the one hand, the equality in~\eqref{exp_chgvar2} is used to obtain
\begin{multline*}
h^2 Q_h^n=\Ka^{n-1/2}u^{P-1}\int_{I_P}\omega_{1-\alpha}(s-x_{n-1})\,ds
\\
+\sum_{j=n+1}^{P}[\Ka^{n+1/2}\fd u^{j}-\Ka^{n-1/2}\fd u^{j-1}]\int_{I_j}\omega_{1-\alpha}(s-x_{n})\,ds,\end{multline*}
where for the second sum, one shows that %for {$\Ka \in C^2(\overline \Omega)$ and $u\in C^3(\overline \Omega)$} it holds
\begin{align*}
\Ka^{n+1/2}\fd u^{j}-\Ka^{n-1/2}\fd &u^{j-1}
= [(\Ka^{n+1/2}-\Ka^n)+\Ka^n][\fd u^{j}-\fd u^{j-1}]+\overline\fd\Ka^n\fd u^{j-1}\\
&=h^2[\frac{h}{2}\Ka'(x_n)+\Ka^n] u''(x_{j-1})+h^2 \Ka'(x_n)u'(x_{j-1})+O(h^3)\\
&=\frac{h^3}{2}\Ka'(x_n)u''(x_{j-1})+h^2\Ka^n [u''(q-h)+(u''(x_{j-1})-u''(q-h))]\\
&\quad +h^2 \Ka'(x_n)[u'(q-h)+(u'(x_{j-1})-u'(q-h))]+O(h^3)\\
&=h^2\Ka^n u''(q-h)+h^2 \Ka'(x_n)u'(q)+O(h^3),
\end{align*}
for any $q \in (x_{j-1},x_j)$.
One the other hand, using equation~\eqref{exp_chgvar1} we have
\[
h\,G_h^n=\sum_{j=n+1}^{P}\int_{I_j}\omega_{1-\alpha}(q-x_{n}) [\Ka^{n}u'(q)-\Ka^{n-1} u'(q-h)]\,dq\\
-\Ka^{n-1}\int_{I_P}\omega_{1-\alpha}(s-x_{n-1})u'(s)\,ds,
\]
where, by Taylor series expansion,
\begin{align*}
\Ka^{n}u'(q)-\Ka^{n-1} u'(q-h)&= \Ka^n[u'(q)-u'(q-h)]+\fd\Ka^{n}\,u'(q-h)\\
&= h\Ka^n u''(q-h)+h\Ka'(x_n)u'(q-h)+O(h^2).
\end{align*}
Inserting the above estimates in \eqref{eq: imp estimate}, we obtain for $1\le n\le P-1$
\[
	T_h^n=E^n+O(h), ~~
	E^n:=-h^{-2}\int_{I_P}\omega_{1-\alpha}(s-x_{n-1})[h\Ka^{n-1}u'(s)+\Ka^{n-1/2}u^{P-1}]ds.
\]
Since %Observing finally that
\[
	\Ka^{n-1/2}u^{P-1}=[\Ka^{n-1}+O(h)][-hu'(x_{P-1})+O(h^2)]=-h \Ka^{n-1}u'(x_{P-1})+O(h^2),
\]
%{for $\Ka\in C^1(\overline I_P)$ and $u\in C^2(\overline I_P)$}, we conclude that
\begin{align*}
	E^n&=O(1)\int_{I_P}{\omega_{1-\alpha}}(s-x_{n-1})ds=O(h)\omega_{1-\alpha}(\xi-x_{n-1}),\quad{\rm for~some}~~\xi \in I_P.
\end{align*}
This completes the proof of the RS truncation error.  $\quad \Box$

\section{Two-sided fractional derivative}\label{sec: two-sided}

In this section, we return to the two-sided fractional differential equation~\eqref{eq: FPDE1}.
To construct our finite difference approximation we simply combine the finite difference schemes introduced in the two previous sections for the LS and RS fractional derivatives.
Specifically, using~\eqref{eq: scheme LS} and~\eqref{eq: scheme RS}, the finite difference solution $U^n\approx u^n$ of the fractional model problem~\eqref{eq: FPDE1} is given by the equations
\begin{multline*}%\label{eq: scheme 2S}
\Ka^{n-1/2}[\theta\,{_aI_x^{1-\alpha}} \partial U(x_{n})+(1-\theta)\,{_xI_b^{1-\alpha}} \partial U(x_{n-1})]\\-
\Ka^{n+1/2}[\theta\,{_aI_x^{1-\alpha}} \partial U(x_{n+1})+(1-\theta)\,{_xI_b^{1-\alpha}} \partial U(x_{n})]=h^2 f^n,
\end{multline*}
for $n=1,\cdots,P-1$, and $U^0=U^P=0.$

The finite difference solution is obtained by solving the linear system ${\bf B}{\bf  U}={\bf F},$ where
${\bf B}=\theta{\bf B}_L+(1-\theta){\bf B}_R$, with the definitions of the matrices ${\bf B}_L$ and ${\bf B}_R$ given in the previous sections.
For instance, for $\theta = 1/2$ we get
\[{\bf B}= \frac{1}{2} \begin{bmatrix}
\ell_{1,1}&\ell_{1,2}&d_{1,3}&d_{1,4}&d_{1,5}&\cdots&d_{1,P-1}\\
\ell_{2,1}&\ell_{2,2}&\ell_{2,3}&d_{2,4}&d_{2,5}&\cdots&d_{2,P-1}\\
c_{3,1}&\ell_{3,2}&\ell_{3,3}&\ell_{3,4}&d_{3,5}&\cdots&d_{3,P-1}\\
c_{4,1}&c_{4,2}&\ell_{4,3}&\ell_{4,4}&\ell_{4,5}&\cdots&b_{4,P-1}\\
\vdots&\vdots&\vdots&\vdots &\vdots&\cdots&\vdots \\
c_{P-1,1}&c_{P-1,2}&c_{P-1,3}&c_{P-1,4}&\cdots&\ell_{P-1,P-2}&\ell_{P-1,P-1}
\end{bmatrix}\]
where $\ell_{i,i}=c_{i,i}+d_{i,i}=(\Ka^{i-1/2}+\Ka^{i+1/2})[3-2^{1-\alpha}]$, and
\begin{align*}
\ell_{i+1,i}&=c_{i+1,i}+d_{i+1,i}=\Ka^{i+1/2}[2^{1-\alpha}-3]-\Ka^{i+3/2}[3^{1-\alpha}-2^{2-\alpha}+1],\\
\ell_{i,i+1}&=c_{i,i+1}+d_{i,i+1}=\Ka^{i+1/2}[2^{1-\alpha}-3]-\Ka^{i-1/2}[3^{1-\alpha}-2^{2-\alpha}+1].
\end{align*}
This shows that the numerical scheme amounts to inverting a system of $(P-1)$ linear equations in the $P-1$ unknowns, so the existence of the finite difference solution follows from its uniqueness. Following a similar path as for the proof of uniqueness for the cases of the LS and RS  fractional derivative schemes (\eqref{eq: matrix1} and \eqref{eq: matrix2}), we obtain
\begin{equation}\label{eq: matrix3}
[\theta {\bf W}_\alpha+ (1-\theta){\bf W}_\alpha^T] \Phi =\psi {\bf K},
\end{equation}
with $\psi=\theta\fd U^1+(1-\theta) \sum_{j=1}^{P}\we_{j,1}\fd U^j.$
%\[
%	\Ka^{n-1/2}\Big(\theta\sum_{j=1}^{n}\we_{n,j}\fd U^j+(1-\theta) \sum_{j=n}^{P}\we_{j,n}\fd U^j\Big)-\Ka^{n+1/2}\Big(\theta\sum_{j=1}^{n+1}\we_{n+1,j}\fd %U^j+(1-\theta)\sum_{j=n+1}^{P}\we_{j,n+1}\fd U^j\Big)=0,\]
%	\[\sum_{n=1}^m \Ka^{n-1/2}\Big(\theta\sum_{j=1}^{n}\we_{n,j}\fd U^j+(1-\theta) \sum_{j=n}^{P}\we_{j,n}\fd %U^j\Big)-\sum_{n=2}^{m+1}\Ka^{n-1/2}\Big(\theta\sum_{j=1}^{n}\we_{n,j}\fd U^j+(1-\theta)\sum_{j=n}^{P}\we_{j,n}\fd U^j\Big)=0,\]
%	\[\Ka^{1/2}\Big(\theta\fd U^1+(1-\theta) \sum_{j=1}^{P}\we_{j,1}\fd U^j\Big)-\Ka^{m+1/2}\Big(\theta\sum_{j=1}^{m+1}\we_{m+1,j}\fd %U^j+(1-\theta)\sum_{j=m+1}^{P}\we_{j,m+1}\fd U^j\Big)=0,\]
By \cite[Lemma A.2]{KimMcLeanMustapha2016}, the  matrix ${\bf W}_\alpha$ in \eqref{eq: Walpha} is positive definite and so is  ${\bf W}_\alpha^T$. Thus, the Toeplitz matrix $\theta{\bf W}_\alpha+(1-\theta){\bf W}_\alpha^T$ is also positive definite and hence, has a inverse, denoted by ${\bf E}_{\alpha,\theta}$, with entries $e_{i,j}$. From \eqref{eq: matrix3}, $\Phi ={\bf E}_{\alpha,\theta}{\bf K}\psi$ and thus, $\sum_{i=1}^P \fd U^i=\psi\sum_{i=1}^P \sum_{j=1}^P e_{i,j}\, k_j.$ Since  $\sum_{j=1}^P\fd U^j=0$,   \begin{equation}\label{eq: SDD1}
\psi\sum_{j=1}^P k_j \sum_{i=1}^P e_{i,j}=0,\quad {\rm where}~~k_j>0.\end{equation}
Recall that, the sequence $\{b_j\}_{j\ge 0}$  is positive, slowly decaying and is also  strictly {\em log-convex}, then by  following the arguments for the case of LS fractional derivative,  we conclude that the matrix ${\bf E}_{\alpha,\theta}$ is strictly diagonally dominant \cite{HornJohnson2013},  $e_{i,i}>0$, and $e_{i,j}\le 0$ for $i\ne j.$ Hence,   $\sum_{i=1}^P e_{i,j}> 0$  and thus,   $\psi=0$ from \eqref{eq: SDD1}. Substitute this  in \eqref{eq: matrix3} yields   $\Phi={\bf 0}$ and it  follows that $U^n=0$ for $1\le n\le P-1$ because $U^0=U^P=0$. This completes the proof of the existence and uniqueness of $U$.

Furthermore, by combining the results of sections~\ref{sec: LS} and~\ref{sec: RS}, it is trivial to show that the truncation error is of order $O(h)$ (not near the boundaries at $x=a,b$), provided that the regularity conditions on $\Ka$, $f$ and $u$ stated in Theorems ~\ref{sec: truncation LS} and~\ref{sec: truncation RS} are met.

 \section{Numerical results}\label{sec: numerical results}
In this section we present several numerical experiments to support the theoretical analyses of the previous sections.
Specifically, we consider the model problem in~\eqref{eq: FPDE1} over $\Omega=(0, 1)$, subject to homogeneous Dirichlet (absorbing) boundary conditions, and  we set $\Ka=1+\exp(x)$.
The finite difference discretization uses uniform spatial meshes with  $P=2^l$ subintervals, for $l>1$, such that $h=1/P$.
The solution error $E_h$ is measured using the discrete $L^\infty$-norm $\|v\|_h=\max_{0\le i\le P}|v(x_i)|$. Based on this error definition, the numerical estimate of convergence rates $\sigma_{h}$ of the finite difference solutions is obtained from the relation $\sigma_h=\log_2(E_{2h}/E_h).$

\begin{table}[hbt]
\caption{Discrete $L^\infty$-norm errors $E_h$ and estimated
numerical convergence rates $\sigma_h$ for different values of $\alpha$, $\theta$ and spatial discretization step size $h$.}
\label{table: errors}
\begin{center}
\renewcommand{\arraystretch}{1.2}
\begin{tabular}{|c|c|cc|cc|cc|}
\hline
& &
\multicolumn{2}{c|}{$\alpha=0.25$}&
\multicolumn{2}{c|}{$\alpha=0.50$}&
\multicolumn{2}{c|}{$\alpha=0.75$}\\
$\theta$&$-\log_2 h$ &     $E_h$ & $\sigma_h$& $E_h$    &$\sigma_h$&  $E_h$    & $\sigma_h$ \\ \hline
\hline
     %       &5&   4.104e-04&         &   3.028e-04&         &   1.771e-04&         \\
            &6&   2.069e-04&   0.9877&   1.568e-04&   0.9493&   9.656e-05&   0.8750\\
            &7&   1.040e-04&   0.9929&   8.028e-05&   0.9659&   5.164e-05&   0.9030\\
       $0.0$&8&   5.214e-05&   0.9960&   4.080e-05&   0.9765&   2.723e-05&   0.9234\\
            &9&   2.611e-05&   0.9976&   2.064e-05&   0.9834&   1.421e-05&   0.9382\\
           &10&   1.307e-05&   0.9986&   1.040e-05&   0.9882&   7.357e-06&   0.9496\\
  \hline\hline
      %       &5&   6.831e-04&         &   3.560e-04&         &   1.488e-04&\\
             &6&   3.528e-04&   0.9535&   1.876e-04&   0.9239&   8.120e-05&   0.8739\\
             &7&   1.784e-04&   0.9838&   9.622e-05&   0.9636&   4.275e-05&   0.9255\\
$0.25$&8       &   8.875e-05&   1.0071&   4.843e-05&   0.9905&   2.200e-05&   0.9588\\
             &9&   4.325e-05&   1.0369&   2.393e-05&   1.0173&   1.108e-05&   0.9887\\
            &10&   2.033e-05&   1.0894&   1.150e-05&   1.0569&   5.432e-06&   1.0290\\
  \hline\hline
  % &5&   9.889e-04&         &   3.774e-04&         &   1.307e-04&\\
   &6&   5.451e-04&   0.8593&   2.024e-04&   0.8990&   7.127e-05&   8.7540\\
   &7&   2.865e-04&   0.9280&   1.045e-04&   0.9530&   3.705e-05&   9.4381\\
$0.5$&8& 1.461e-04&   0.9713&   5.269e-05&   0.9883&   1.868e-05&   9.8776\\
   &9&   7.289e-05&   1.0036&   2.599e-05&   1.0198&   9.157e-06&   1.0287\\
  &10&   3.545e-05&   1.0398&   1.243e-05&   1.0643&   4.304e-06&   1.0893\\
\hline\hline
 %  &5&   6.378e-04&         &   3.409e-04&         &  1.427e-04&\\
   &6&   3.353e-04&   0.9282&   1.818e-04&   0.9071&  7.899e-05&   0.8529\\
   &7&   1.714e-04&   0.9672&   9.392e-05&   0.9527&  4.190e-05&   0.9147\\
$0.75$&8&8.632e-05&   0.9898&   4.757e-05&   0.9812&  2.167e-05&   0.9513\\
   &9&   4.289e-05&   1.0092&   2.370e-05&   1.0054&  1.097e-05&   0.9820\\
  &10&   2.094e-05&   1.0341&   1.156e-05&   1.0359&  5.408e-06&   1.0205\\
  \hline\hline
 %  &5&   4.018e-04&         &   2.926e-04&         &   1.687e-04&\\
   &6&   2.047e-04&   0.9728&   1.537e-04&   0.9289&   9.350e-05&   0.8512\\
   &7&   1.034e-04&   0.9855&   7.929e-05&   0.9546&   5.048e-05&   0.8893\\
$1.0$&8& 5.197e-05&   0.9922&   4.048e-05&   0.9700&   2.677e-05&   0.9149\\
   &9&   2.607e-05&   0.9956&   2.053e-05&   0.9794&   1.403e-05&   0.9326\\
  &10&   1.306e-05&   0.9975&   1.037e-05&   0.9857&   7.283e-06&   0.9457\\
  \hline
     \end{tabular}
\end{center}
\end{table}

{\bf Example 1.} We first consider the source term $f$ leading to the exact solution
\begin{equation}\label{eq: exact}
u_{\rm ex}(x)=x^{4-\theta(1-\alpha)}(1-x)^{4-(1-\theta)(1-\alpha)}.
\end{equation}
%The determination of  $f$ corresponding to $u_{\rm ex}$ is detailed in Appendix ~\ref{appen: RHS}.

\begin{figure}[hbt]
\begin{center}
\includegraphics[scale=0.4]{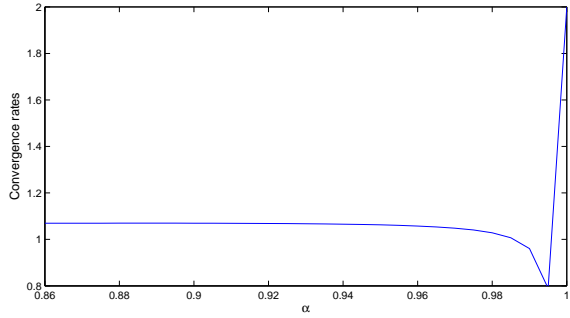}
\end{center}
\caption{Graphical plot of the numerical convergence rates $\sigma_{h}$ against the diffusion exponent $\alpha$. Computations use $\theta=1/2$ and $P=1024.$}
\label{fig: 1}
\end{figure}

We first fix $\theta=1/2$, $P=1024$ and report in Fig.~\ref{fig: 1} the estimates $\sigma_h$ as a function of $\alpha$. The plot shows that $\sigma_h\sim 1$, denoting an error in $O(h)$, for almost all values of $\alpha$ except in the immediate neighborhood of $\alpha=1$. When $\alpha\rightarrow 1$, $\sigma_h$ exhibits a rapidly varying behavior to reach the expected second order convergence rate  at $\alpha=1$.

\begin{figure}[hbt]
\begin{center}
\includegraphics[scale=0.5]{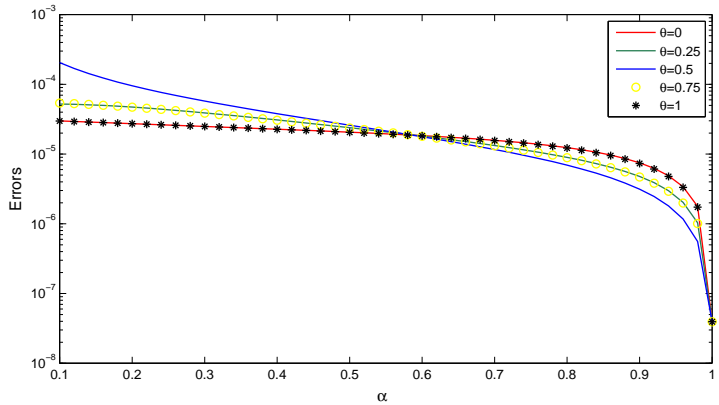}
\end{center}
\caption{The  error $E_h$ against the diffusion exponent $\alpha$, for  $P = 512$ and different values of $\theta$ as indicated.}
\label{fig: 2}
\end{figure}

Next, we fix $P=512$ and plot  $E_h$ against $\alpha$ for different values of $\theta$. Results are reported in Fig.~\ref{fig: 2}.
We observe that the errors are almost the same for $\theta=0.25$ and $\theta=0.75$, and for $\theta=0$ and $\theta=1$. This is due to the similar singularity behavior near the boundaries of the exact solution  in~\eqref{eq: exact} for any choice of $\theta =c$ and $\theta =1-c$. Note that the errors are decreasing as $\alpha \rightarrow 1$ for all $\theta$.
Interestingly enough, Fig.~\ref{fig: 2} also shows that for $\alpha < 0.6$, the error is lower for extreme values of $\theta$, that is close to 0 or 1, and on the contrary $E_h$ is lower for intermediate values ($\approx 1/2$) when $\alpha > 0.6$.

Table~\ref{table: errors} reports the $L^\infty$-norm of $E_h$ and the corresponding estimates of convergence rate for different values of $\alpha$, $\theta$ and the discretization step size $h$. The table confirms the $O(h)$ errors, for all the values of $\alpha$ and $\theta$ shown, as $h$ goes to zero.

{\bf Example 2.}  (non-smooth solutions) In practice, due to the presence of the two-sided fractional derivative, the solution $u$ of  \eqref{eq: FPDE1} admits end-point singularities  even if the source term $f$ is  smooth.
It was proved recently in~\cite{MaoChenShen2016} that, for $\theta=1/2$, the leading singularity term takes the form $x^{\frac{1+\alpha}{2}}(1-x)^{\frac{1+\alpha}{2}}$ when the diffusivity coefficient $\Ka$ is constant.
Similarly, one can show that  leading singularity term takes the form $(x-a)^{\alpha}$, with $a=0$ presently, in the case of LS fractional derivative ($\theta=1$), and the form $(b-x)^{\alpha}$, with $b=1$ presently, in the case of  RS fractional derivatives ($\theta=0$).
For smooth $\Ka$, we conjecture the same singular behavior.
Furthermore, we suggest that for $\theta \in [0,1]$,  the leading singularity term has the generic form $(x-a)^{1-\theta(1-\alpha)}(b-x)^{1-(1-\theta)(1-\alpha)}$ ($a=0$ and $b=1$). However, demonstrating this point remains an open problem and it will be a subject of future work. Noting that, for $\kappa=1$ and for $0\le \theta \le 1$,  the authors in \cite{ErvinHeuerRoop2017} studied the regularity properties of the solution $u$ of problem \eqref{eq: FPDE1} where the fractional derivative operator is not of the Riemann-Liouville type, see \cite[Equations (1.3) and (3.11)]{ErvinHeuerRoop2017}.

To support our claim, we choose now the source term $f$ such that $u_{\rm ex}(x)=x^{1-\theta(1-\alpha)}(1-x)^{1-(1-\theta)(1-\alpha)}$ is the exact solution of the problem with other settings as before.
One can easily check that the truncation errors analyses provided above are not valid in this situation.
We then apply to this problem our finite difference scheme for the LS ($\theta=1$) and RS ($\theta=0$) fractional derivatives cases for different values of $\alpha$ and $h$. Table~\ref{table: nonsmooth} reports the discrete $L^\infty$-norm of $E_h$ and estimates of the convergence rates $\sigma_h$. The results clearly indicate a convergence rate of the error in $O(h^\alpha)$.

\begin{table}[hbt]
\caption{Discrete $L^\infty$-norm errors $E_h$ and estimated
numerical convergence rates $\sigma_h$ for different values of $\alpha$, $\theta$ and spatial discretization step size $h$.}
\label{table: nonsmooth}
\begin{center}
\renewcommand{\arraystretch}{1.2}
\begin{tabular}{|c|c|cc|cc|cc|} \hline
& &
\multicolumn{2}{c|}{$\alpha=0.25$}&
\multicolumn{2}{c|}{$\alpha=0.50$}&
\multicolumn{2}{c|}{$\alpha=0.75$}\\
$\theta$&$-\log_2 h$ &     $E_h$ & $\sigma_h$& $E_h$    &$\sigma_h$&  $E_h$    & $\sigma_h$ \\ \hline
\hline
   %   &6&  6.120e-02&        &  2.733e-02&         &  7.845e-03&        \\
   %   &7&  5.057e-02&  0.2752&  1.916e-02&   0.5123&  4.624e-03&  0.7626\\
      &8&  4.214e-02&  0.2632&  1.348e-02&   0.5068&  2.732e-03&  0.7590\\
$0.0$ &9&  3.527e-02&  0.2567&  9.510e-03&   0.5037&  1.618e-03&  0.7556\\
     &10&  2.959e-02&  0.2534&  6.716e-03&   0.5019&  9.601e-04&  0.7533\\
     &11&  2.485e-02&  0.2517&  4.745e-03&   0.5010&  5.702e-04&  0.7518\\
     &12&  2.088e-02&  0.2509&  3.354e-03&   0.5005&  3.388e-04&  0.7510\\
  \hline \hline
  %     &6&   5.739e-02&         &   2.638e-02&         &   7.452e-03&         \\
  %     &7&   4.895e-02&   0.2297&   1.881e-02&   0.4877&   4.496e-03&   0.7289\\
       &8&   4.145e-02&   0.2399&   1.336e-02&   0.4940&   2.692e-03&   0.7402\\
$1.0$  &9&   3.498e-02&   0.2449&   9.465e-03&   0.4970&   1.606e-03&   0.7454\\
      &10&   2.946e-02&   0.2475&   6.700e-03&   0.4985&   9.562e-04&   0.7478\\
      &11&   2.480e-02&   0.2487&   4.740e-03&   0.4993&   5.690e-04&   0.7490\\
      &12&   2.086e-02&   0.2494&   3.352e-03&   0.4996&   3.384e-04&   0.7495\\
      \hline
\end{tabular}
\end{center}
\end{table}

\begin{figure}[hbt]
\begin{center}
\includegraphics[scale=0.5]{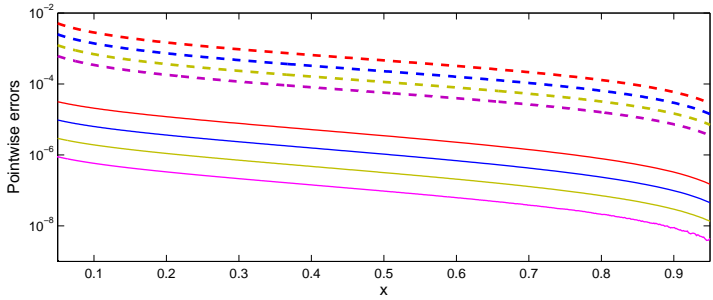}
\end{center}
\caption{Pointwise errors using uniform (dashed lines) and nonuniform meshes with $\gamma=3$  (solid
lines), for α = 0.25 with $P = 256, 512, 1024, 2048$ (in order from top to bottom).}
\label{fig: 3}
\end{figure}

This degradation of the convergence rate was expected because the low regularity of the solution: $u_{\rm ex}\in C^\alpha[0,1]$.
In the context of time-stepping schemes for fractional diffusion of fractional wave equations, adapted meshes with refinement (clustering of elements) around the singularity successfully improve the errors and consequently, the convergence rates, see~\cite{McLeanMustapha2007,Mustapha2011}.
To check if such refinement approach could be useful in our problem of (steady) spatial fractional diffusion problem, we set $\theta=1$ (LS singularity) and consider a family of graded spatial meshes of $\Omega=(0,1)$ based on a sequence of points given by $x_i = (i/P)^\gamma$, $i=0,\dots,P$ and $\gamma\ge 1$ is a refinement parameter.
The objective is to refine the mesh at the boundary $x=0$ where the solution has a singularity. Table~\ref{table: graded meshes} reports the evolution with $\log_2(P)$ of the $L^\infty$-norm of the error and estimated convergence rate $\sigma_h$ and using $\gamma=2$, 3 and 4. The results show that one can obtain a convergence rate of the error that is $O(h^{\alpha \gamma})$.
Finally, Fig.~\ref{fig: 3} compares the pointwise errors obtained for uniform and non-uniform meshes with $\gamma=3$ when using the same number of discretization points $P=256$, 512, 1024 and 2048.
The reduction of the error due to the mesh refinement is clearly visible.
Note that similar results can be obtained for $\theta=0$ using discretization points defined by $x_i=1-((P-i)/P)^\gamma$ to refine the mesh at the endpoint $x=1$.

\begin{table}[hbt]
\caption{Discrete $L^\infty$-norm errors $E_h$ and estimated
numerical convergence rates $\sigma_h$ for $\alpha=0.25$, $\theta=1$ (LS fractional derivatives), different number of discretization points ($P$) and refinement parameters $\gamma$.}
\label{table: graded meshes}
\begin{center}
\renewcommand{\arraystretch}{1.2}
\begin{tabular}{|c|cc|cc|cc|} \hline
&
\multicolumn{2}{c|}{$\gamma=2$}&
\multicolumn{2}{c|}{$\gamma=3$}&
\multicolumn{2}{c|}{$\gamma=4$}\\
$\log_2 P$ &     $E_h$ & $\sigma_h$& $E_h$    &$\sigma_h$&  $E_h$    & $\sigma_h$ \\ \hline
\hline
   %    6&   2.300e-02&         &   8.128e-03&         &   2.871e-03&        \\
   %    7&   1.628e-02&   0.4988&   4.838e-03&   0.7484&   1.438e-03&   0.9976\\
       8&   1.151e-02&   0.4996&   2.878e-03&   0.7495&   7.194e-04&   0.9992\\
       9&   8.140e-03&   0.4998&   1.711e-03&   0.7498&   3.597e-04&   0.9997\\
      10&   5.756e-03&   0.4999&   1.018e-03&   0.7499&   1.800e-04&   0.9999\\
      11&   4.070e-03&   0.4999&   6.051e-04&   0.7499&   8.994e-05&   0.9999\\
      12&   2.878e-03&   0.5002&   3.600e-04&   0.7500&   4.497e-05&   0.9998\\
      \hline
    \end{tabular}
\end{center}
\end{table}

%%%%%%%%%%%%%%%%%%%%%%%%%%%%%%%%%%%%%%%%%%%%%%%%%%%%%%%%%%%%%%%%%%%%%%%%%%%%%%%%%%%%%%%%%%%%%%%%%%%%%%%%%%5
 \section{Concluding remarks}\label{sec: conclusion}
The objective of this work was to propose and analyze a finite-difference scheme for the solution of general one-dimensional fractional elliptic problems with a variable diffusion coefficient. For the proposed scheme, we proved the existence and uniqueness of the numerical  solution and established the order of convergence for the truncation error with the spatial step size. Some numerical results were also presented for problems admitting both smooth and nonsmooth solutions.

This paper will form a stepping stone for the researchers who are interested in computational solutions of variable coefficient two-sided fractional derivative problems.
The results obtained in this work lead to several questions that will have to be addressed in the future.
First, it will crucial to address the reason(s) for the dramatic deterioration in the order of convergence of the finite difference scheme when the fractional order $\alpha$ immediately departs from $1$ (classical case)?
Second, it will be interesting to explore the possibility of incorporating the fractional exponent $\alpha$ directly in the finite difference discretization, that is, fractionalizing the numerical scheme. A possible route along this direction could be inspired by the recent research papers on the fractionalization of the Crank-Nicolson time-scheme for solving time-fractional diffusion equation, see~\cite{Dimitrov2014}.
Finally, mechanisms for determining the order of singularity near the boundaries in the case of variable diffusivity remains to be developed. A possibility could be to look at series solution to~\eqref{eq: FPDE1}.  These and other related open questions will be the subject of future research.

\end{document}